\documentclass[journal]{IAENGtran}
\usepackage{cite}
\ifCLASSINFOpdf
   \usepackage[pdftex]{graphicx}
   \usepackage{epstopdf}%% or use the graphicx package for more complicated commands
   \usepackage{epsfig}%% or use the epsfig package if you prefer to use the old commands
   \usepackage{amsthm}%% The amsthm package provides extended theorem environments
   \usepackage{amsmath}
   \usepackage{multirow}
   \usepackage[utopia]{mathdesign}
   \usepackage{fancyhdr}
   \usepackage[bookmarksnumbered,bookmarksopen,colorlinks=true,citecolor=black,linkcolor=black]{hyperref}
   \DeclareGraphicsExtensions{.pdf,.jpeg,.png}
\else
   \usepackage[dvips]{graphicx}
   \DeclareGraphicsExtensions{.eps}
\fi
\usepackage[font=footnotesize]{subfig}
\usepackage{stfloats}
\theoremstyle{plain}

\theoremstyle{plain}

\theoremstyle{plain}
\newtheorem{algorithm}{Algorithm}

\theoremstyle{remark}

\begin{document}
\title{The nonlinear HSS-like iterative method for absolute value equations}
%
%
% author names and IAENG memberships
\author{Mu-Zheng Zhu$^{\dagger}$ ~\IAENGmembership{Member,~IAENG,}
        and~Ya-E Qi$^{\ddagger}$ %~\IAENGmembership{Senior Member,~IAENG,}
%        and~Guo-Feng Zhang$^{*}$%~\IAENGmembership{Fellow,~IAENG}% <-this % stops a space
\thanks{Manuscript received January 2, 2018;\  % revised June XX, 20XX.
%(Write the date on which you submitted your paper for review.)
This work was supported by the National Natural Science Foundation of China(11661033) and the Scientific Research Foundation for Doctor of Hexi University.}
%(sponsor and financial support acknowledgment goes here).
%Papertitles should be written in uppercase and lowercase letters, not all uppercase. Avoid writing long formulas with subscripts in the title; short formulas that identify the elements are fine.}
\thanks{$^{\dagger}$ Mu-Zheng Zhu is with the School of Mathematics and Statistics, Lanzhou University, Lanzhou, 734000 P. R. China. And he is also with the School of Mathematics and Statistics, Hexi University, Zhangye, 734000 P. R. China;  e-mail: zhumzh07@yahoo.com.}% <-this % stops a space
\thanks{$^{\ddagger}$ Ya-E Qi is with the School of Chemistry and Chemical Engineering, Hexi University, Zhangye 734000 P.R. China.}
%\thanks{$^{*}$ Guo-Feng Zhang is corresponding author, with School of Mathematics and Statistics, Lanzhou University, Lanzhou, 730000 P. R. China; e-mail: gf\_zhang@lzu.edu.cn.}
}% <-this % stops a space

\maketitle

\pagestyle{empty}
\thispagestyle{empty}

\begin{abstract}
%\boldmath
Salkuyeh proposed the Picard-HSS iteration method to solve the absolute value equation (AVE), which is a class of non-differentiable NP-hard problem. To further improve its performance, a nonlinear HSS-like iteration method is proposed. Compared to that the Picard-HSS method is an inner-outer double-layer iteration scheme, the HSS-like iteration is only a monolayer and the iteration vector could be updated timely. Some numerical experiments are used to demonstrate that the nonlinear HSS-like method is feasible, robust and effective.
\end{abstract}
% IAENGtran.cls defaults to using nonbold math in the Abstract.
% of the abstract to achieve this. Many IAENG journals frown on math

\begin{IAENGkeywords}
absolute value equation, nonlinear HSS-like iteration, fixed point iteration, positive definite.
\end{IAENGkeywords}

% \fi
%
% For peerreview papers, this IAENGtran command inserts a page break and
\IAENGpeerreviewmaketitle

\section{Introduction}\label{sec:1}
% \IAENGPARstart{A}{demo} file is ....
% normal text (unknown if ever used by IAENG):
% \IAENGPARstart{A}{}demo file is ....
% \IAENGPARstart{T}{his demo} file is ....
\IAENGPARstart{T}{HE} solution of the absolute value equation (AVE) of the following form is considered:
\begin{equation}\label{eq:1}
 Ax-|x|=b.
\end{equation}
Here, $A\in \mathbb {C}^{n\times n}$, $x,\,b\in \mathbb {C}^{n}$ and $|x|$ denotes the component-wise absolute value of vector $x$, i.e., $|x|=(|x_1|,\,|x_2|,\,..., |x_n|)^T$.
The AVE (\ref{eq:1}) is a special case of the generalized absolute value equation (GAVE) of the type
\begin{equation}\label{eq:1a}
 Ax-B|x|=b,
\end{equation}
where $A, B\in \mathbb {C}^{m\times n}$ and $x,\ b\in \mathbb {C}^{m}$. The GAVE (\ref{eq:1a}) was introduced in \cite{rohn2004theorem} and investigated in a more general context in \cite{hu2010note,mangasarian2012primal,mangasarian2006absolute}. Recently, these problems have been investigated in the literature \cite{mangasarian2007absolute,mangasarian2006absolute,noor2012iterative,prokopyev2009equivalent,rohn2014iterative}.

The AVE (\ref{eq:1}) arises in linear programs, quadratic programs, bimatrix games and other problems, which can all be reduced to a linear complementarity problem (LCP) \cite{cottle2009linear,mangasarian1977solution}, and the LCP is equivalent to the AVE (\ref{eq:1}). This implies that AVE is NP-hard in its general form \cite{mangasarian2007absolute,mangasarian2006absolute,noor2012iterative}.
Beside, if $B=0$, then the generalized AVE (\ref{eq:1a}) reduces to a system of
linear equations $Ax = b$, which have many applications in scientific computation \cite{noor2012iterative}.

The main research of AVE includes two aspects: one is the theoretical analysis, which focuses on the theorem of alternatives, various equivalent reformulations, and the existence and nonexistence of solutions; see \cite{caccetta2011globally,hu2010note,prokopyev2009equivalent,rohn2004theorem}.
And the other is how to solve the AVE. We mainly pay attention to the letter.

In the last decade, based on the fact that the LCP is equivalent to the AVE and the special structure of AVE, a large variety of methods for solving AVE (\ref{eq:1}) can be found in the literature; See \cite{mangasarian2009knapsack,mangasarian2012primal,noor2012iterative,rohn2014iterative}.  These also include the following: a finite succession of linear programs (SLP) is established in \cite{mangasarian2007absolute,mangasarian2009knapsack}, which arise from a reformulation of the AVE as the minimization of a piecewise-linear concave function on a polyhedral set and solving the latter by successive linearization; a semi-smooth Newton method is proposed, which largely shortens the computation time than the SLP method in \cite{mangasarian2009generalized}; %Then, to formulate the NP-hard n-dimensional knapsack feasibility problem as an equivalent absolute value equation in an n-dimensional noninteger real variable space, Mangasarian further proposed a finite succession of linear programs (SLP) for solving the AVE in \cite{mangasarian2009knapsack}.
Furthermore, a smoothing Newton algorithm is presented in \cite{caccetta2011globally}, which is proved to be globally convergent and the convergence rate is quadratic under the condition that the singular values of $A$ exceed 1. This condition is weaker than the one used in \cite{mangasarian2009generalized}.
%a generalized Newton method, which has global and finite convergence, was proposed for the AVE in \cite{zhang2009global}. This method utilizes both the semismooth and the smoothing Newton steps, in which the semismooth Newton step guarantees the finite convergence and the smoothing Newton step contributes to the global convergence. And

Recently, The Picard-HSS iteration method is proposed to solve AVE by Salkuyeh in \cite{salkuyeh2014picard}, which is originally designed to solve weakly nonlinear systems \cite{bai2009hss} and its generalizations are also paid attention \cite{zhu2013modified,pu2013class}. The sufficient conditions to guarantee the convergence of this method and some numerical experiments are given to show the effectiveness of the method. However, the numbers of the inner HSS iteration steps are often problem-dependent and difficult to be determined in actual computations. Moreover, the iteration vector can not be updated timely. In this paper, we present the nonlinear HSS-like iteration method to overcome the defect mentioned above, which is designed originally for solving weakly nonlinear systems in \cite{bai2009hss}.

The rest of this paper is organized as follows. In Section \ref{sec:2} the HSS and Picard-HSS iteration methods are reviewed. In Section \ref{sec:3} the nonlinear HSS-like iteration method for solving AVE (\ref{eq:1}) is described. Numerical experiments are presented in Section \ref{sec:4}, to shown the feasibility and effectiveness of the nonlinear HSS-like method.  Finally, some conclusions and an open problem are drew in Section \ref{sec:5}.

\section{The HSS and Picard-HSS iteration methods}\label{sec:2}
\IAENGPARstart{I}{N} this section, the HSS iteration method for solving the non-Hermitian linear systems and the Picard-HSS iteration method for solving the AVE (\ref{eq:1}) are reviewed.
%\subsection{The HSS iteration method}

Let $A\in \mathbb {C}^{n \times n}$ be a non-Hermitian positive definite matrix, $B\in \mathbb {C}^{n \times n}$ be a zero matrix, the GAVE (\ref{eq:1a}) reduced to the non-Hermitian system of linear equations
\begin{equation}\label{eq:5}
Ax=b.
\end{equation}

Because any square matrix $A$ possesses a Hermitian and skew-Hermitian splitting (HSS)
\begin{equation}\label{eq:6}
A=H+S,\quad H=\frac{1}{2}(A+A^H)\quad \text{and}\quad S=\frac{1}{2}(A-A^H),
\end{equation}
the following HSS iteration method is first introduced by Bai, Golub and Ng in \cite{bai2003hermitian} for the solution of the non-Hermitian positive definite system of linear equations (\ref{eq:5}).

\begin{algorithm} (The HSS iteration method.)\\
Given an initial guess $x^{(0)}\in \mathbb {C}^{n}$, compute $x^{(k)}$ for $k =0,\ 1,\ 2,\ ...$ using the following iteration scheme until $\{x^{(k)}\}_{k =0}^{\infty }$ converges,
\begin{equation*}
\begin{aligned}
\left\{
\begin{array}{ll}
 (\alpha I+H)x^{(k +\frac{1}{2})}=(\alpha I-S)x^{(k )}+b, \\
 (\alpha I+S)x^{(k +1)}=(\alpha I-H)x^{(k +\frac{1}{2})}+b,
\end{array}
\right.
\end{aligned}
\end{equation*}
where $\alpha $ is a positive constant and $I$ is the identity matrix.
\end{algorithm}

%In the matrix-vector form, the above HSS iteration can be equivalently rewritten as
%\begin{equation*}
%x^{(k +1)}= T(\alpha )x^{(k )}+G(\alpha )b = T(\alpha )^{k +1} x^{(0)}+\sum _{j=0}^{k } T(\alpha )^j G(\alpha ) b, \quad k =0,1,2,... ,
%\end{equation*}
%where
%\begin{equation*}
%T(\alpha )=(\alpha I+S)^{-1}(\alpha I-H)(\alpha I+H)^{-1} (\alpha I-S)
%\quad\text{and}\quad
%G(\alpha )=2\alpha (\alpha I+S)^{-1}(\alpha I+H)^{-1}.
%\end{equation*}
%Here, $T(\alpha )$ is the iteration matrix of the HSS method.

When the matrix $A\in\mathbb{C}^{n\times n}$ is positive definite, i.e. its Hermitian part $H=\frac{1}{2}(A+A^H)$ is positive definite, Bai et al. proved that the spectral radius of the HSS iteration matrix is less than 1 for any positive parameters $\alpha$, i.e., the HSS iteration method is unconditionally convergent; see  \cite{bai2003hermitian}.
%
%\subsection{The Picard-HSS iteration methods}

For the convenience of the subsequent discussion, the AVE (\ref{eq:1}) can be rewritten as its equivalent form:
\begin{equation*}
Ax=f(x), \quad f(x)=|x|+b.
\end{equation*}

Recalling that the linear term $Ax$ and the nonlinear term $f(x)=|x|+b$ are well separated and the Picard iteration method is a fixed-point iteration, the Picard iteration
%\begin{equation*}
$Ax^{(k+1)}=f(x^{(k)}),\quad k=0,\ 1,\ ....,$
%\end{equation*}
can be used to solve the AVE (\ref{eq:1}).
When the matrix $A\in\mathbb{C}^{n\times n}$ is large sparse and positive definite, the next iteration $x^{(k+1)}$ may be inexactly computed by HSS iteration. This naturally lead to the following iteration method proposed in \cite{salkuyeh2014picard} for solving the AVE (\ref{eq:1}).

\begin{algorithm} (The Picard-HSS iteration method)\\
Let $A\in \mathbb {C}^{n \times n}$ be a sparse and positive definite matrix, $H=\frac{1}{2}(A+A^H)$ and $S=\frac{1}{2}(A-A^H)$ be its Hermitian and skew-Hermitian parts respectively. Given an initial guess $x^{(0)}\in \mathbb {C}^n$ and a sequence $\{\ell_k\}_{k=0}^{\infty }$ of positive integers, compute $x^{(k+1)}$ for $k=0,1,2,\ldots $ using the following iteration scheme until $\{x^{(k)}\}$ satisfies the stopping criterion:

(a) Set $x^{(k,0)}:=x^{(k)};$

(b) For $\ell =0,1,\ldots ,\ell_k-1$, solve the following linear systems to obtain $x^{(k,\ell +1)}$:
\begin{equation*}
\begin{aligned}
\left\{
\begin{array}{ll}
(\alpha I+H)x^{(k,\,\ell +\frac{1}{2})}=(\alpha I-S)x^{(k,\,\ell )}+|x^{(k)}|+b,\\
(\alpha I+S)x^{(k,\,\ell +1)}=(\alpha I-H)x^{(k,\,\ell +\frac{1}{2})}+|x^{(k)}|+b,
\end{array}
\right.
\end{aligned}
\end{equation*}
where $\alpha $ is a given positive constant and $I$ is the identity matrix;

 (c) Set $x^{(k+1)}:=x^{(k,\ell_k)}$.
\end{algorithm}

The advantage of the Picard-HSS iteration method is obvious. First, the two linear sub-systems in all inner HSS iterations have the same shifted Hermitian coefficient matrix $\alpha I +H$ and shifted skew-Hermitian coefficient matrix $\alpha I+S$, which are constant with respect to the iteration index $k$. Second, As the coefficient matrix $\alpha I +H$ and $\alpha I+S$ are Hermitian and skew-Hermitian respectively, the first sub-system can be solved exactly by making use of the Cholesky factorization and the second one by the LU factorization. The last, these two sub-systems can be solve approximately by the conjugate gradient method and a Krylov subspace method like GMRES, respectively; see \cite{bai2009hss,salkuyeh2014picard}.

\section{The nonlinear HSS-like iteration method}\label{sec:3}

\IAENGPARstart{I}{N} the Picard-HSS iteration, the numbers $\ell_k,\ k=0,\ 1, 2,\ ...$ of the inner HSS iteration steps are often problem-dependent and difficult to be determined in actual computations \cite{bai2009hss}. Moreover, the iteration vector can not be updated timely.
Thus, to avoid these defect and still preserve the advantages of the Picard-HSS iteration method, based on the HSS (\ref{eq:6}) and the nonlinear fixed-point equations
\begin{equation*}
(\alpha I+H)x=(\alpha I-S)x+|x|+b,
\end{equation*}
and
\begin{equation*}
(\alpha I+S)x=(\alpha I-H)x+|x|+b,
\end{equation*}
the following nonlinear HSS-like iteration method is proposed to solve the AVE (\ref{eq:1}).

\begin{algorithm} {The nonlinear HSS-like iteration method.}\\
Let $A\in \mathbb {C}^{n \times n}$ be a sparse and positive definite matrix, $H=\frac{1}{2}(A+A^H)$ and $S=\frac{1}{2}(A-A^H)$ be its Hermitian and skew-Hermitian parts respectively. Given an initial guess $x^{(0)}\in \mathbb {C}^n$, compute $x^{(k+1)}$ for $k=0,1,2,\ldots $ using the following iteration scheme until $\{x^{(k)}\}$ satisfies the stopping criterion:
\begin{equation*}\label{eq:hsslike}
\begin{aligned}
\left\{
\begin{array}{ll}
(\alpha I+H)x^{(k+\frac{1}{2})}=(\alpha I-S)x^{(k)}+|x^{(k)}|+b,\\
(\alpha I+S)x^{(k+1)}=(\alpha I-H)x^{(k+\frac{1}{2})}+|x^{(k+\frac{1}{2})}|+b,
\end{array}
\right.
\end{aligned}
\end{equation*}
where $\alpha $ is a given positive constant and $I$ is the identity matrix.
\end{algorithm}

It is obvious that both $x$ and $|x|$ in the second step are updated in the nonlinear HSS-like iteration, but only $x$ is updated in the Picard-HSS iteration. Furthermore, the nonlinear HSS-like iteration is a monolayer iteration scheme, but the Picard-HSS is an inner-outer double-layer iteration scheme.

To obtain a one-step form of the nonlinear HSS-like iteration, we define
\begin{equation*}
\begin{aligned}
%\left\{
\begin{array}{ll}
U(x)&=(\alpha I +H)^{-1}((\alpha I-S)x+|x|+b),\\
V(x)&=(\alpha I +S)^{-1}((\alpha I-H)x+|x|+b),
\end{array}
%\right.
\end{aligned}
\end{equation*}
and
\begin{equation*}
\psi(x):=V\circ U(x)=V(U(x)).
\end{equation*}
Then the nonlinear HSS-like iteration scheme can be equivalently expressed as
\begin{equation*}\label{eq:iter1}
x^{(k+1)}=\psi(x^{(k)}).
\end{equation*}

The Ostrowski theorem, i.e., Theorem 10.1.3 in \cite{ortega2000iterative}, gives a local convergence theory about a one-step stationary nonlinear iteration. Based on this, Bai et al. established the local convergence theory for the nonlinear HSS-like iteration method in \cite{bai2009hss}. However, these convergence theory  has a strict requirement that $f(x)=|x|+b$ must be $F$-differentiable at a point $x^*\in\mathbb{D}$ such that $Ax^*-|x^*|=b$. Obviously, the absolute value function $|x|$ is non-differentiable.
Thus, the convergence analysis of the nonlinear HSS-like iteration method for solving weakly nonlinear linear systems is unsuitable for solving AVE, and need further discuss.

At the end of this section, we remark that the main steps in the nonlinear HSS-like iteration method can be alternatively reformulated into residual-updating form as follows.

\begin{algorithm} (The HSS-like iteration method (residual-updating variant).) \ Given an initial guess $x^{(0)}\in \mathbb{D}\subset\mathbb {C}^n$, compute $x^{(k+1)}$ for $k=0,1,2,\ldots $ using the following iterative procedure until $\{x^{(k)}\}$ satisfies the stopping criterion:

(1) Set: $r^{(k)}:=|x^{(k)}|+b-Ax^{(k)}$,

(2) Solve: $(\alpha I+H)v=r^{(k)}$,

(3) Set: $x^{(k+\frac{1}{2})}=x^{(k)}+v$,\ $r^{(k)}:=|x^{(k+\frac{1}{2})}|+b-Ax^{(k+\frac{1}{2})}$,

(4) Solve: $(\alpha I+S)v=r^{(k)}$,

(5) Set: $x^{(k+1)}=x^{(k+\frac{1}{2})}+v$,\\
where $\alpha $ is a given positive constant and $I$ is the identity matrix.
\end{algorithm}

\section{Numerical experiments}\label{sec:4}

\IAENGPARstart{I}{N} this section, the numerical properties of the Picard, Picard-HSS and nonlinear HSS-like methods are examined and compared experimentally by a suit of test problems. All the tests are performed in MATLAB R2013a on Intel(R) Core(TM) i5-3470 CPU 3.20 GHz and 8.00 GB of RAM, with machine precision $10^{-16}$, and terminated when the current residual satisfies
\begin{equation*}
%\text{RES}=
\frac{\Vert Ax^{(k)}-|x^{(k)}|-b \Vert _{2}}{\Vert b\Vert _{2}}\leq 10^{-5},
\end{equation*}
where %RES is the norm of absolute residual vectors,
$x^{(k)}$ is the computed solution by each of the methods at iteration $k$, and a maximum number of the iterations 500 is used.

In addition, the stopping criterion for the inner iterations of the Picard-HSS method is set to be
\begin{equation*}
\frac{\Vert b^{(k)}-A\,s^{(k,\,\ell_k )}\Vert _2}{\Vert b^{(k)}\Vert _2}\leq\eta_k,
\end{equation*}
where $b^{(k)}=|x^{(k)}|+b-Ax^{(k)}$, $s^{(k,\,\ell_k )}=x^{(k,\,\ell_k)}-x^{(k,\,\ell_k-1)}$, $\ell_k$ is the number of the inner iteration steps and $\eta_k$ is the prescribed tolerance for controlling the accuracy of the inner iterations at the $k$-th outer iteration. If $\eta_k$ is fixed for all $k$, then it is simply denoted by $\eta$. Here, we take $\eta=0.1$.
%In our numerical experiments, the accuracy of the inner iterations $\eta_k$ for the Picard-HSS iteration method is fixed and set to $\eta=0.1$.

The first subsystem with the Hermitian positive definite coefficient matrix $(\alpha I+H)$ in (\ref{eq:hsslike}) is solved by the Cholesky factorization, and the second subsystem with the skew-Hermitian coefficient matrix $(\alpha I+S)$ in (\ref{eq:hsslike}) is solved by the LU factorization.

The optimal parameters employed in the Picard-HSS and nonlinear HSS-like iteration methods have been obtained experimentally. In fact, the experimentally found optimal parameters are the ones resulting in the least numbers of iterations and CPU times\cite{salkuyeh2014picard}. As mentioned in \cite{bai2009hss} the computation of the optimal parameter is often problem-dependent and generally difficult to be determined.

We consider the two-dimensional convection-diffusion equation

\begin{equation*}
\begin{cases}
-(u_{xx}+u_{yy})+q(u_x+u_y)+p u =f(x,y), \  (x,y)\in \Omega , \\
u(x,y)=0, \  (x,y)\in \partial \Omega ,
\end{cases}
\end{equation*}
where $\Omega =(0,1)\times (0,1)$, $\partial \Omega$ is its boundary, $q$ is a positive constant used to measure the magnitude of the diffusive term and $p$ is a real number. We use the five-point finite difference scheme to the diffusive terms and the central difference scheme to the convective terms. Let $h=1/(m+1)$ and $Re=(qh)/2$ denote the equidistant step size and the mesh Reynolds number, respectively. Then we get a system of linear equations $Ax=d$, where $A$ is a matrix of order $n=m^2$ of the form

\begin{equation}\label{eq:ex}
A=T_x \otimes I_m+I_m \otimes T_y+p I_n,
\end{equation}
with
$$T_x=\mathrm{tridiag}(t_2,\,t_1,\,t_3)_{m\times m} \  \mathrm{and} \ T_{y}=\mathrm{tridiag}(t_2,\,0,\,t_3)_{m\times m},$$
where $t_1=4,\ t_2=-1-Re,\ t_3=-1+Re$, $I_m$ and $I_n$ are the identity matrices of order $m$ and $n$ respectively, $\otimes $ means the Kronecker product.
%, and $T_x$ and $T_y$ are the tridiagonal matrices
%
%\begin{equation}\label{eq:ex1}
%T_x=\mathrm{tridiag}(t_2,\,t_1,\,t_3)_{m\times m}\quad \mathrm{and} \quad T_{y}=\mathrm{tridiag}(t_2,\,0,\,t_3)_{m\times m} ,
%\end{equation}
%with
%
%\begin{equation}\label{eq:ex2}
%t_1=4,\quad t_2=-1-Re,\quad t_3=-1+Re.
%\end{equation}

In our numerical experiments, the matrix $A$ in AVE (\ref{eq:1}) is defined by (\ref{eq:ex}) with different values of $q\,(q=0,\,1,\,10,\,100\, \text{and}\,1000)$ and different values of $p\,(p=0\ \text{and}\ 0.5)$. It is easy to find that for every nonnegative number $q$ the matrix $A$ is in general non-symmetric positive definite\cite{salkuyeh2014picard}. We use the zero vector as the initial guess, and the right-hand side vector $b$ of AVE (\ref{eq:1}) is taken in such a way that the vector $x=(x_1,x_2,\ldots ,x_n)^T$ with $x_k=(-1)^k\ {\rm i}\ (k=1,\,2,\,\ldots,\,n)$ is the exact solution, where $\rm{i}$ denotes the imaginary unit.

\begin{table}[!hbp]
\tabcolsep 3mm
\renewcommand{\arraystretch}{1.2}
\caption{The optimal parameters values $\alpha$\ (p=0).}\label{tab:1}
\begin{tabular}{lllrrrrrrrrrr}
\hline
\multicolumn{2}{c}{Optimal parameters}  &m=10    &m=20    &m=40    &m=80    \\
\hline
q=0 &HSS-like  &1.3& 1.0& 1.0& 1.0\\
    &Picard-HSS&1.1& 0.5& 0.2& 0.1\\
q=1 &HSS-like  &1.4& 1.0& 1.0& 1.0 \\
    &Picard-HSS&1.1& 0.6& 0.3& 0.2\\
q=10&HSS-like  &1.7& 1.1& 1.0& 1.0 \\
   &Picard-HSS &1.6& 0.8& 0.4& 0.2\\
q=100&HSS-like &2.5& 2.7& 1.7& 1.2\\
    &Picard-HSS&2.4& 2.7& 1.8& 0.9\\
%q=1000&HSS-like&1.9& 1.1& 2.9& 2.3 \\
%   &Picard-HSS&1.9& 1.1& 2.9& 2.3\\
\hline
\end{tabular}
\end{table}

\begin{table}[!hbp]
\tabcolsep 3mm
\renewcommand{\arraystretch}{1.2}
\caption{The optimal parameters values $\alpha$ (p=0.5).}\label{tab:2}
\begin{tabular}{lllrrrrrrrrrr}
\hline
\multicolumn{2}{c}{Optimal parameters}  &m=10    &m=20    &m=40    &m=80    \\
\hline
q=0 &HSS-like  &2.4 & 2.2 & 2.1 & 2.0\\
    &Picard-HSS&2.2 & 2.0 & 1.8 & 1.8\\
q=1 &HSS-like  &2.4 & 2.2 & 2.1 & 2.0 \\
    &Picard-HSS&2.3 & 2.0 & 1.8 & 1.8\\
q=10&HSS-like  &2.6 & 2.3 & 2.2 & 2.1 \\
   &Picard-HSS &2.4 & 2.3 & 2.0 & 1.9\\
q=100&HSS-like &3.4 & 2.9 & 2.3 & 2.3\\
    &Picard-HSS&3.5 & 3.0 & 2.3 & 2.1\\
%q=1000&HSS-like&2.9 & 2.4 & 2.4 & 2.5\\
%   &Picard-HSS &2.8 & 2.4 & 2.4 & 2.5\\
\hline
\end{tabular}
\end{table}

\begin{table}[!hbp]
\renewcommand{\arraystretch}{1.2}
\caption{Numerical results for test problems with different values of $m$ and $q$ ($p=0$, RES($\times 10^{-6}$)\ ).}\label{tab:3}
\begin{tabular}{lllrrrrrrrrrr}
\hline
    &Methods   &   &m=10    &m=20    &m=40    &m=80    \\
\hline
q=0 &HSS-like  &IT &27    &35    &65    &81    \\
    &          &CPU&0.0375    &0.0146    &0.1016    &0.6085    \\
    &          &RES&9.4084    &8.7487    &9.9395    &9.9502    \\
    &Picard-HSS&IT$_\text{out}$&5    &5    &5    &5    \\
    &          &IT$_\text{int}$&7.2  &13.8 &33   &62.6 \\
    &          &IT &36    &69    &165    &313    \\
    &          &CPU&0.0084    &0.0250    &0.2310    &2.0708    \\
    &          &RES&5.2907    &7.1401    &7.9627    &9.1458    \\
    &Picard    &IT &--    &--    &--    &--    \\
    &          &CPU&--    &--    &--    &--    \\
    &          &RES&--    &--    &--    &--    \\
\hline
q=1 &HSS-like  &IT &28    &38    &65    &81    \\
    &          &CPU&0.0044    &0.0199    &0.1343    &0.8436    \\
    &          &RES&8.7445    &9.5272    &9.9148    &9.9588    \\
    &Picard-HSS&IT$_\text{out}$&5    &5    &5    &5    \\
    &          &IT$_\text{int}$&7.2  &13.6 &27   &64.8  \\
    &          &IT &36    &68    &135    &324    \\
    &          &CPU&0.0050    &0.0317    &0.2527    &3.0404    \\
    &          &RES&6.3073    &8.0703    &7.7121    &9.3360    \\
    &Picard    &IT &--    &--    &--    &--    \\
    &          &CPU&--    &--    &--    &--    \\
    &          &RES&--    &--    &--    &--    \\
\hline
q=10&HSS-like  &IT &17    &32    &51    &85    \\
    &          &CPU&0.0029    &0.0176    &0.1077    &0.8857    \\
    &          &RES&7.8979    &7.2166    &9.3825    &9.8324    \\
    &Picard-HSS&IT$_\text{out}$&5    &5    &5    &5    \\
    &          &IT$_\text{int}$&3.8  &7    &13.2  &25.4 \\
    &          &IT &19    &35    &66    &127    \\
    &          &CPU&0.0031    &0.0174    &0.1285    &1.2305    \\
    &          &RES&2.6888    &4.0994    &5.9529    &7.1369    \\
    &Picard    &IT &--    &--    &--    &--    \\
    &          &CPU&--    &--    &--    &--    \\
    &          &RES&--    &--    &--    &--    \\
\hline
q=100&HSS-like &IT &18    &20    &25    &42    \\
    &          &CPU&0.0037    &0.0117    &0.0574    &0.4687    \\
    &          &RES&8.2690    &8.8682    &7.5469    &9.3710    \\
    &Picard-HSS&IT$_\text{out}$&5    &5    &5    &5    \\
    &          &IT$_\text{int}$&3.8  &4.2  &5.6  &8.2 \\
    &          &IT &19    &21    &28    &41    \\
    &          &CPU&0.0039    &0.0116    &0.0602    &0.4413    \\
    &          &RES&3.2385    &3.3042    &3.3640    &3.6666    \\
    &Picard    &IT &4    &8    &39    &--    \\
    &          &CPU&0.0009    &0.0036    &0.0436    &--    \\
    &          &RES&6.9831    &0.0032    &6.8249    &--    \\
%q=1000&HSS-like&IT &21    &38    &50    &51    \\
%    &          &CPU&0.0045    &0.0330    &0.2909    &2.3601    \\
%    &          &RES&9.3457    &8.6466    &8.4651    &8.9837    \\
%    &Picard-HSS&IT$_\text{out}$&5    &5    &5    &5    \\
%    &          &IT$_\text{int}$&5    &7.6  &10   &10\\
%    &          &IT &25    &38    &50    &50    \\
%    &          &CPU&0.0049    &0.0318    &0.2872    &2.3036    \\
%    &          &RES&0.8263    &3.4996    &6.6572    &8.5153    \\
%    &Picard    &IT &3    &4    &5    &12    \\
%    &          &CPU&0.0009    &0.0047    &0.0442    &0.7562    \\
%    &          &RES&1.4196    &0.0006    &0.0338    &0.0005    \\
\hline
\end{tabular}
\end{table}

\begin{table}[!hbp]
\renewcommand{\arraystretch}{1.2}
\caption{Numerical results for test problems with different values of $m$ and $q$ ($p=0.5$, RES($\times10^{-6}$)\ ).}\label{tab:4}
\begin{tabular}{lllrrrrrrrrrr}
\hline
    &Methods    &   &m=10    &m=20    &m=40    &m=80    \\
\hline
q=0	&HSS-like	&IT &29	&38	&36	&35	\\
	&	        &CPU&0.0037 	&0.0155 	&0.0590 	&0.2849 	\\
	&	        &RES&7.7828 	&8.0756 	&9.6565 	&8.8724 	\\
	&Picard-HSS	&IT$_\text{out}$&5	&5	&5	&5	\\
	&	        &IT$_\text{int}$&7  &14.6&35&66.4\\
    &           &IT             &35	&73	&175	&332	\\
	&	        &CPU&0.0040 	&0.0261 	&0.2420 	&2.2039 	\\
	&	        &RES&5.4444 	&7.4483 	&8.1466 	&9.3423 	\\
	&Picard	    &IT &9	&--	&--	&--	\\
	&	        &CPU&0.0010	&--	&--	&--	\\
	&	        &RES&0.0016	&--	&--	&--	\\
\hline
q=1	&HSS-like	&IT &29	&42	&38	&36	\\
	&	        &CPU&0.0044 	&0.0218 	&0.0824 	&0.4113 	\\
	&	        &RES&8.1442 	&8.5129 	&9.8553 	&8.2976 	\\
	&Picard-HSS	&IT$_\text{out}$&5	&5	&5	&5	\\
	&	        &IT$_\text{int}$&7.8&14.4&28&42 \\
    &           &IT             &39	&72	&140	&210	\\
	&	        &CPU&0.0052 	&0.0335 	&0.2612	&1.9946	\\
	&	        &RES&4.2330 	&5.3548 	&8.8367	&8.5786\\
	&Picard	    &IT &9	&--	&--	&--	\\
	&	        &CPU&0.0011	&--	&--	&--	\\
	&	        &RES&0.0011	&--	&--	&--	\\
\hline
q=10&HSS-like	&IT &18	&34	&45	&42	\\
	&	        &CPU&0.0030 	&0.0183 	&0.0960 	&0.4848 	\\
	&	        &RES&8.1728 	&6.0961 	&8.8821 	&9.2731 	\\
	&Picard-HSS	&IT$_\text{out}$&5	&5	&5	&5	\\
	&	        &IT$_\text{int}$&4  &7  &13.6&25 \\
    &           &IT             & 20&35	&68	&125	\\
	&	        &CPU&0.0032 	&0.0179 	&0.1379 	&1.2244	\\
	&	        &RES&1.5905 	&5.9853 	&5.1449 	&8.9996	\\
	&Picard	    &IT &7	&--	&--	&--	\\
	&	        &CPU&0.0009 	&--	&--	&--	\\
	&	        &RES&0.1525 	&--	&--	&--	\\
\hline
q=100&HSS-like	&IT &14	&14	&22	&37	\\
	&	        &CPU&0.0032 	&0.0088 	&0.0518 	&0.4204 	\\
	&	        &RES&9.8625 	&5.9430 	&5.6508 	&7.4515 	\\
	&Picard-HSS	&IT$_\text{out}$&5	&5	&5	&5	\\
	&	        &IT$_\text{int}$&3.4&3.2&5.5&8.4\\
    &           &IT             & 17&16	&22	&42	\\
	&	        &CPU&0.0037 	&0.0093 	&0.0498 	&0.4558 	\\
	&	        &RES&1.4643 	&1.2321 	&2.7830 	&4.4858 	\\
	&Picard	    &IT &4	&6	&14	&--	\\
	&	        &CPU&0.0009 	&0.0030 	&0.0205 	&--	\\
	&	        &RES&0.9480 	&0.0162 	&0.0229 	&--	\\
\hline
%q=1000&HSS-like	&IT &21	&38	&19	&18	\\
%	&	        &CPU&0.0044 	&0.0330 	&0.1312 	&1.1059 	\\
%	&	        &RES&7.9169 	&8.2226 	&8.5076 	&6.6481 	\\
%	&Picard-HSS	&IT$_\text{out}$&5	&5	&5	&5	\\
%	&	        &IT$_\text{int}$&4.6&7.6&4.4&4  \\
%    &           &IT             &23	&38	&22	&20	\\
%	&	        &CPU&0.0047 	&0.0322 	&0.1452 	&1.1742 	\\
%	&	        &RES&2.2895 	&5.6813 	&1.4223 	&2.0359 	\\
%	&Picard	    &IT &3	&4	&5	&7	\\
%	&	        &CPU&0.0009 	&0.0046 	&0.0419 	&0.6358 	\\
%	&	        &RES&0.7292 	&0.0002 	&0.0003 	&0.3046 	\\
\end{tabular}
\end{table}

In Tables \ref{tab:3} and \ref{tab:4}, we present the numerical results with respect to the Picard, Picard-HSS and nonlinear HSS-like iterations, the experimentally optimal parameters used in the Picard-HSS and nonlinear HSS-like iterations are those given in Tables \ref{tab:1} and \ref{tab:2}. we give the elapsed CPU time in seconds for the convergence (denoted as CPU), the norm of absolute residual vectors (denoted as RES), and the number of outer, inner and total iteration steps (outer and inner iterations only for Picard-HSS) for the convergence (denoted as IT$_{\text{out}}$, IT$_{\text{int}}$ and IT, respectively).
%The number of inner iteration steps for Picard-HSS is a average value of all the inner iterative step corresponded to outer iterative step, and
The number of outer iteration steps for Picard-HSS and the number of iteration steps for Picard and HSS-like iteration methods larger than 500 are simply listed by the symbol "--".

From these two tables, we see that both the HSS-like and Picard-HSS methods can successfully produced approximate solution to the AVE for all of the problem-scales $n=m^2$ and the convective measurements $q$, while the Picard iteration converges only for some special cases. Here, it is necessary to mention that the shifted matrices $\alpha I+H$ and $\alpha I+S$ are usually more well-conditioned than the matrix $A$\cite{salkuyeh2014picard}.

For the convergent cases, the number of iteration steps for the Picard and HSS-like methods and the number of inner iteration steps for the Picard-HSS method are increase rapidly with the increasing of problem-scale, while the number of outer iteration steps is fixed. The CPU time also increases rapidly with the increasing of the problem-scale for all iteration methods.

When the convective measurements $q$ become large, for all iteration method, both the number of iteration steps (except outer iteration for Picard-HSS) and the amount of CPU times decrease, while $q=1000$ is in the opposite situation.

Clearly, in terms of iteration step, the nonlinear HSS-like method and the Picard-HSS are more robust than Picard, and the nonlinear HSS-like method performs much better than the Picard-HSS;
In terms of CPU time, the situation is almost the same, but the Picard iteration method is the most time-efficient in the convergent cases, e.g. $q=1000$.
Therefore, the nonlinear HSS-like method are the winners for solving this test problem when the convective measurements $q$ is small.

%As mentioned above, we have observed that Picard iteration does not successfully compute any solution for small $q$. A possible illustration for this numerical phenomenon is that for large $q$ the matrix property of $A$ may be comparable, or the latter is even more dominant than the former.

\section{Conclusions}\label{sec:5}
\IAENGPARstart{I}{N} this paper we have studied the nonlinear HSS-like iteration method for solving the absolute value equation (AVE). This method is based on separable property of the linear term $Ax$ and nonlinear term $|x|+b$ and the Hermitian and skew-Hermitian splitting of the involved matrix $A$. Compared to that the Picard-HSS iteration scheme is an inner-outer double-layer iteration scheme, the nonlinear HSS-like iteration is a monolayer and the iteration vector could be updated timely. %By leveraging the smoothing approximate function, the locally convergence have been analysed. Further
Numerical experiments have shown that the nonlinear HSS-like method is feasible, robust and efficient nonlinear solver. The most important is it can outperform the Picard-HSS in actual implementation.

\section*{Acknowledgements}
The author would like to thank the anonymous referees for his/her careful reading of the manuscript and useful comments and improvements.

\bibliographystyle{IAENGtran}
\end{document}